\documentclass[11pt,a4paper]{article}

\usepackage{amsmath}
\usepackage{amsfonts}
\usepackage[english]{babel}
\usepackage[latin1]{inputenc}
\usepackage{fontenc}
\usepackage{euscript}
\usepackage{vmargin}
\usepackage{amssymb}
\usepackage{latexsym}
\usepackage{enumerate}
\usepackage[dvips]{graphicx}
\usepackage{hyperref}
\def\rain{\to +\infty}

\def\N{{\rm I\kern-.20em N}}
\def\R{{\rm I\kern-.20em R}}
\def\indi{{1\kern-.20em\rm I}}

\def\bkR{{\rm I\kern-.17em R}}
\def\bkN{{\rm I\kern-.20em N}}

\newtheorem{lema}{Lemma}[section]
\newtheorem{teo}{Theorem}[section]
\newtheorem{nota}{Remark}
\newtheorem{coro}{Corollary}[section]
\newtheorem{prop}{Proposition}[section]

\newtheorem{defi}{Definition}[section]

\newcommand{\pg}{\hspace{0.6cm}}

\newcommand{\bdem} {\begin{proof}}
\newcommand {\edem}{\hfill $\square$ \end {proof}}

\begin{document}
\title{Estimating the Upcrossings Index}
%\author{
\author{J.R. Sebasti\~ao \footnote{ Escola Superior de Gest\~ao,
Instituto Polit\'ecnico de Castelo Branco, 6000 Castelo Branco,
Portugal.\newline
E-mail:jrenato@esg.ipcb.pt}\\ Escola Superior de Gest\~ao\\
Instituto Polit\'ecnico de Castelo Branco\\ Portugal \and A.P.
Martins, H. Ferreira, L. Pereira  \footnote{Departamento de
Matem\'atica, Universidade da Beira Interior, 6200 Covilhã,
Portugal.
E-mail:amartins@mat.ubi.pt; lpereira@mat.ubi.pt; helena@mat.ubi.pt}\\
Departamento de Matem\'atica \\
Universidade da Beira Interior\\ Portugal}
\date{}
\maketitle

\noindent {\bf Abstract:}  For stationary sequences,  under general local and asymptotic dependence restrictions, any limiting point process for time normalized upcrossings of high levels is a compound Poisson process, i.e., there is a clustering of high upcrossings, where the underlying Poisson points represent cluster positions, and the multiplicities correspond to cluster sizes. For such classes of stationary sequences there exists the upcrossings index $\eta,$ $0\leq \eta\leq 1,$ which is directly related to the extremal index $\theta,$ $0\leq \theta\leq 1,$ for suitable high levels. In this paper we consider the problem of estimating the upcrossings index $\eta$ for a class of stationary sequences satisfying a mild oscillation restriction. For the proposed estimator, properties such as consistency and asymptotic normality are studied. Finally, the performance of the estimator is assessed through  simulation studies for autoregressive processes and case studies in the fields of environment and finance.

%In the present paper we introduce an estimator  of the upcrossings index, $\eta,$ which is a measure of clustering of
%upcrossings of a high level by the variables of a stationary
%sequence. The reciprocal of the upcrossing index can be interpreted
%as the limiting mean size of a cluster of upcrossings or as the
%limiting mean size of a cluster of exceedance runs. Identifyiny
%these clusters is a key issue to estimate $\eta.$ If the stationary
%sequence satisfies a local dependence condition, that restricts the
%possibility of having upcrossings separated by non-upcrossings
%within a given block,  there is always a group of successive
%upcrossings in the neighborhood of an upcrossing. Hence, in this
%case, it is natural to estimate $\eta$ as the ratio between the
%total number of non-upcrossing followed by an upcrossing and the
%total number of upcrossings. It is shown that such an estimate is
%typically weakly consistent and asymptotically normal. Finally, the
%results of a simulation study pertaining to this estimator are
%presented for two stationary autoregressive processes.
\vspace{0.3cm}


\vspace{0.3cm}


%% --------------------------------------------------------------
%% Inicio do texto do artigo
%% --------------------------------------------------------------

\section{Introduction and preliminary results}\setcounter{equation}{0}

\pg Let ${\bf{X}}=\{X_n\}_{n\geq 1}$ be a stationary sequence and
 ${\bf{u}}=\{ u_{n}\} _{n\geq 1}$ a sequence of real levels. We say that
  $\bf{X}$ has an upcrossing of $u_{n}$, at $i$, if the event
  $\{X_i\leq u_n<X_{i+1}\},$ $i=1,\ldots,n,$
occurs. The point process of upcrossings of $u_n$ by the first $n$
variables of ${\bf{X}}$ is then defined by
\begin{equation}
\widetilde{N}_n(B)=\sum_{i=1}^{n-1} \indi_{\{X_i\leq u_n<X_{i+1}\}}\delta_{\frac{i}{n}}(B),\qquad B
\subset [0,1],\ n\geq 1,\label{p_up}
\end{equation}
 where $\delta_a(\cdot)$ denotes unit mass at
$a\in B$ measure and  $\indi_A$ the indicator of event  $A$.

Ferreira \cite{fer1} showed that if ${\bf{X}}$ satisfies the mixing $\Delta({\bf{u}})$
condition, introduced in Hsing {\it{et al.}} \cite{hsing}, and
$\widetilde{N}_n$ converges in distribution (as a point process on
[0,1]), then the limit is necessarily a compound Poisson process.
For independent and identically distributed (i.i.d.) sequences it is easy to see that $\widetilde{N}_n$
converges in distribution if and only if $u_n\equiv u_n^{(\nu)}$ for some
$\nu >0,$ where $u_n^{(\nu)}$ are normalized levels for upcrossings,
that is, ${\bf{u}}^{(\nu)}=\{u_n^{(\nu)}\}_{n\geq 1}$ denotes a
sequence that satisfies
\begin{equation}
\lim_{n\to +\infty}nP(X_1\leq u_n^{(\nu)}<X_2)=\nu,\label{eq1}
\end{equation}
being the limit point process a Poisson process on [0,1] with intensity
$\nu.$ We remark that when ${\bf{X}}$ is a stationary sequence, satisfying
the long range dependence condition  $D({\bf{u}}^{(\nu)})$ of
Leadbetter \cite{lead1} and the mild oscillation restriction
$D''({\bf{u}}^{(\nu)})$ of Leadbetter and Nandagopalan \cite{lead2},
$\widetilde{N}_n$ converges in distribution to the same Poisson
process as in the i.i.d. case (Leadbetter and Nandagopalan \cite{lead2}). When ${\bf{X}}$ only satisfies
condition $\Delta({\bf{u}}^{(\nu)})$ the limit compound Poisson process
$\widetilde{N}$ is in fact a direct consequence of the clustering of
high level upcrossings in a dependent sequence. Its intensity is simply
$\eta \nu,$ where $\eta$ is the upcrossings index, and has
multiplicity distribution $\widetilde{\pi}(j)=\lim_{n\to
+\infty}\widetilde{\pi}_n(j),$ $j=1,2,\ldots,$ with
\begin{equation}
\widetilde{\pi}_n(j)=P\left(\sum_{i=1}^{r_n}\indi_{\{X_i\leq u_n<X_{i+1}\}}=j\
\Bigl|\ \sum_{i=1}^{r_n}\indi_{\{X_i\leq u_n<X_{i+1}\}}>0\right),\
j=1,2,\ldots\label{pi}
\end{equation} for some sequence $r_n=[n/k_n]$ and $\{k_n\}_{n\geq
1}$ a sequence of positive integers satisfying
\begin{equation}
k_n\xrightarrow [n\rain]{}+\infty ,\qquad \frac{k_nl_n}{n}
\xrightarrow [n\rain]{}0, \qquad k_n\alpha_{n,l_n}\xrightarrow
[n\rain]{} 0,\label{eq2}
\end{equation}
where $\alpha_{n,l_n}$ are the mixing coefficients of the $\Delta({\bf{u}})$ condition.

The upcrossings of $u_n$ by $X_i$ with $i\in
J_{n,j}=\{(j-1)r_n+1,\ldots,jr_n\},$ for some $j=1,\ldots,k_n,$ are
regarded as forming a cluster of upcrossings and $\widetilde{\pi}_n(j),\
j=1,2,\ldots$ is called the conditional cluster size distribution, since it is simply the distribution of the number of upcrossings in a cluster (i.e. in a set $J_{n,j}$) given that there is at least one.

The upcrossings index can then be viewed as a measure of clustering of upcrossings of high levels $u_n$ by the variables in $\bf{X}$ and is formally defined as follows.

\begin{defi} [Ferreira \cite{fer1}] If for each $\nu>0$ there exists ${\bf{u}}^{(\nu)}=\{u_n^{(\nu)}\}_{n\geq 1}$ and\linebreak
$\displaystyle{\lim_{n\to +\infty}
P(\widetilde{N}^{(\nu)}_n([0,1])=0)=e^{-\eta \nu}},$ for
some constant $0\leq \eta \leq 1,$ then we say that the sequence
${\bf{X}}$ has upcrossings index $\eta.$
\end{defi}

Many common cases, as i.i.d. sequences or sequences that satisfy condition $D''({\bf{u}}^{(\nu)})$, have upcrossings index $\eta=1.$ A value $\eta<1$ indicates clustering of upcrossings of ${\bf{u}}^{(\nu)},$ giving rise to multiplicities in the limit.

If for each $\nu>0$ there exists $u_n^{(\nu)}$ satisfying (\ref{eq1}) and $\lim_{n\to +\infty}P(X_1>u_n^{(\nu)})=\tau$ for some $\tau>0$ then the upcrossings index $\eta$ exists if and only if the
extremal index $\theta$ exists  and, in this case,
\begin{equation}
\theta=\frac{\nu}{\tau}\eta.\label{rel_eta_teta}
\end{equation} (Ferreira \cite{fer1}). Note that the previous conditions imply that $\theta\leq \eta$  since if a level $u_n$ is simultaneously normalized for upcrossings and for exceedances we necessarily have $\nu\leq \tau$.

The conditional cluster size distribution $\widetilde{\pi}_n$ defined in (\ref{pi}) is also related to the upcrossings index $\eta$ in the following  interesting manner, which is a direct consequence of Lemma 2.1 of Ferreira \cite{fer1}.

\begin{prop}
Suppose ${\bf{X}}$ satisfies $\Delta({\bf{u}}^{(\nu)}),$ $\nu>0,$ has upcrossings index $\eta>0$ and let $\widetilde{\pi}_n$ be defined as in (\ref{pi}) with $u_n=u_n^{(\nu)}.$ Then $$\lim_{n\to +\infty}\sum_{j\geq 1} j\widetilde{\pi}_n(j)=\frac{1}{\eta}.$$
\end{prop}

The reciprocal of the upcrossings index can then be interpreted as the limiting mean cluster size of
upcrossings. Hence, it is natural to estimate $\eta$ as the reciprocal of
the sample average cluster size. However, such an estimator suffers
from the drawback that the identification of clusters (equivalently,
choice of $k_n$) depends on the knowledge of the mixing coefficients
$\alpha_{n,l_n}.$ Alternative ways of identifying clusters of high level upcrossings is therefore a key issue for the estimation of $\eta.$

If a sequence is suitably well-behaved then one might hope that
groups of successive upcrossings of $u_n$ are sufficiently far apart
that each group can be regarded as a separate cluster. A sufficient
condition for this to hold is condition $\widetilde{D}^{(3)}({\bf{u}})$ of Ferreira \cite{fer1}.
%\begin{equation}
%\lim_{n\to +\infty}n \sum_{j=4}^{r_n}P(X_1\leq u_n<X_2,\
%\widetilde{N}_{3,3}=0,\ X_j\leq u_n<X_{j+1}) =0.\label{cond}
%\end{equation}

\begin{defi} [Ferreira \cite{fer1}]
Let ${\bf{X}}$ be a sequence satisfying the condition
$\Delta({\bf{u}}).$ ${\bf{X}}$ satisfies
condition $\widetilde{D}^{(3)}({\bf{u}})$ if $$\lim_{n\to
+\infty}nP(X_1\leq u_n<X_2,\ \widetilde{N}_{3,3}=0,\
\widetilde{N}_{4,r_n}>0)=0,$$ for some sequence $r_n=[n/k_n]$ with
$\{k_n\}_{n\geq 1}$ satisfying (\ref{eq2}) and $\widetilde{N}_{i,j}\equiv\widetilde{N}_n([i/n,j/n])$ with $\widetilde{N}_{i,j}=0$ if $j<i.$
\end{defi}

\begin{nota}
Condition $\widetilde{D}^{(3)}({\bf{u}})$ is implied by condition
$$
\lim_{n\to \infty}n\sum_{j=4}^{r_n-1}P(X_1\leq u_n<X_2, \widetilde{N}_{3,3}=0,X_{j}\leq u_n<X_{j+1})=0,
$$ which is clearly implied by condition
\begin{equation}
\lim_{n\to \infty}n\sum_{j=4}^{r_n-1}P(X_1\leq u_n<X_2,,X_{j}\leq u_n<X_{j+1})=0. \label{condU}
\end{equation}
\end{nota}

\begin{nota}
Sequences that satisfy condition $D''({\bf{u}})$ also satisfy condition $\widetilde{D}^{(3)}({\bf{u}})$ and the latter belongs to a wide family of local dependence conditions $\widetilde{D}^{(k)}({\bf{u}}),\ k\geq 2,$ introduced in Ferreira \cite{fer1}. This family of local dependence conditions is slightly stronger than the family of conditions $D^{(k)}({\bf{u}}),\ k\geq 1,$ considered in Chernick {\it{et al.}} \cite{cher}. The negatively correlated uniform AR(1) process given in Chernick et al. \cite{cher} is an example of a process that verifies the previous condition as shown in Sebastião {\it{et al.}} \cite{seb}. In Ferreira \cite{fer1} we find an example of a max-autoregressive process for which $\widetilde{D}^{(3)}({\bf{u}})$ also holds.
\end{nota}

\begin{nota}
Necessary and sufficient conditions for (\ref{condU}) to hold, involving upcrossings-tail dependence coefficients, can be found in Ferreira and Ferreira \cite{fer3}.
\end{nota}
%This condition implies condition $\widetilde{D}^{(3)}({\bf{u}})$ of
%Ferreira (2006)

%Under condition $\widetilde{D}^{(3)}({\bf{u}}),$ $\lim_{n\to +\infty}P(\widetilde{N}([0,1])=0)$ can
%be computed from the joint distribution of 4 consecutive terms of
%$X$.

%\begin{prop}[Ferreira (2006)]
%Suppose that conditions $\Delta({\bf{u}})$ and  $\widetilde{D}^{(3)}({\bf{u}})$ hold
%for ${\bf{X}}$ and \linebreak $\displaystyle{\lim_{n\to
%+\infty}P(\widetilde{N}([0,1])=0)}>0.$ Then
%$$P(\widetilde{N}([0,1])=0)-\exp(-nP(X_1\leq u_n<X_2,\ \widetilde{N}_{3,3}=0))\xrightarrow
%[n\rain]{} 0.$$
%\end{prop}

Condition $\widetilde{D}^{(3)}({\bf{u}})$  locally restricts the dependence of the sequence, but still allows clustering of upcrossings. It roughly states that whenever an upcrossing of a high level occurs, a cluster of upcrossings may follow it, but once the sequence falls below the threshold it is very unlikely to upcross it again in the nearby observations. Thus, it enables the identification of upcrossings clusters since for this class upcrossings may be simply identified asymptotically as runs of consecutive upcrossings and the cluster sizes as run lengths. Indeed, Ferreira \cite{fer2} proved that if the  conditional upcrossing run length distribution is defined as
\begin{equation}
\widetilde{\pi}_n^*(j)=P(X_3\leq u_n<X_4,\ldots,X_{2j+1}\leq u_n
<X_{2j+2},\widetilde{N}_{2j+3,2j+3}=0\ |\
\widetilde{N}_{1,1}=0,X_3\leq u_n<X_4),\ j\geq 1,\label{run_prob}
\end{equation}
then the following result gives the conditional expected length of an upcrossing run.

\begin{prop}[Ferreira \cite{fer2}]
If, for each $\nu >0,$ ${\bf{X}}$ satisfies condition
$\widetilde{D}^{(3)}({\bf{u}}^{(\nu)})$  and $P(X_3\leq u_n^{(\nu)}<X_4,\ldots,X_{2j-1}\leq u_n^{(\nu)}
<X_{2j})\longrightarrow 0,$ as $j\to +\infty,$ then the upcrossings
index of ${\bf{X}}$ exists and is equal to $\eta$ if and only if,
\begin{equation}
\sum_{j\geq 1} j\widetilde{\pi}_n^*(j)=\frac{P(X_1\leq u_n^{(\nu)}<X_2)}{P(X_1\leq
u_n^{(\nu)}<X_2,\widetilde{N}_{3,3}=0)}\xrightarrow
[n\rain]{}
\frac{1}{\eta},\label{eq3}
\end{equation}
for each $\nu>0.$
\end{prop}

From this result it seems natural to estimate $\eta$ by the reciprocal of the sample average run length. We shall consider such an estimator in Section 3.
%Whereas, for stationary time series, the extremal index $\theta$ summarizes short-term extremal behaviour during a stress period, known as volatility clustering, the upcrossings index can be viewed as a measure of duration of high volatility near extremes, with $\eta$ decreasing with the increase of duration of high volatility.
%Attending to this result the upcrossings index can be in some sense interpreted as the reciprocal of the ``mean time of duration of high volatility between extreme events'' and it is directly related to the upcrossings of high levels.

\begin{nota}
Under condition $\widetilde{D}^{(3)}({\bf{u}})$  the upcrossings index can be related to other dependence measures, namely the lag-1 upcrossings tail dependence coefficient $\mu_1,$ since $\eta=1-\mu_1=1-\lim_{x\to x_F} P(X_3\leq u_n^{(\nu)}<X_4\ |\ X_1\leq u_n^{(\nu)}<X_2 ),$ where $x_F$ is the upper limit of the common marginal distribution $F$ of ${\bf{X}},$ as shown in Ferreira and Ferreira \cite{fer3}.
\end{nota}

Alternative expressions for $\eta$ involving only stationarity are given in the following simple lemma.

\begin{lema}
If $nP(X_1\leq u_n<X_2,\widetilde{N}_{3,3}=0)\xrightarrow[n\rain]{} \xi>0$ then the following are equivalent:
\begin{enumerate}
\item[\bf{i)}] $P(\widetilde{N}_{3,3}=0\ |\ X_1\leq u_n<X_2)\xrightarrow[n\rain]{} \eta;$
\item[\bf{ii)}] $nP(X_1\leq u_n<X_2)\xrightarrow[n\rain]{} \frac{\xi}{\eta};$ %(i.e. $u_n^{(\nu)}=u_n^{(\xi/\eta)});$
\item[\bf{iii)}] $n(1-P(\widetilde{N}_{1,1}=0,\ \widetilde{N}_{3,3}=0))\xrightarrow[n\rain]{} \xi+\frac{\xi}{\eta}.$
\end{enumerate}
\end{lema}

\begin{nota}
From the stationarity of ${\bf{X}}$ we obviously have $P(X_1\leq u_n<X_2,\ \widetilde{N}_{3,3}=0)=P(\widetilde{N}_{1,1}=0,\ X_3\leq u_n<X_4)$ and consequently an upcrossings run can be either identified at its beginning or at its end.
\end{nota}

The upcrossings index estimation is important not only by itself, as a measure of clustering of upcrossings of high levels, but also because of its relation with other dependence coefficients, namely the extremal index and the lag-1 upcrossings tail coefficient (Ferreira and Ferreira \cite{fer3}).
%if we assume that a sea-wall projection requires a coastal defense from all sea-levels then the estimation of the probability of upcrossings of the sea-wall height my be of interest.
Therefore, we shall pursuit this issue in the remainder of this paper. Note that under the conditions of Proposition 1.2 it is reasonable to estimate $\eta$ by the ratio between the number of upcrossings followed by a non-upcrossing and the number of upcrossings of a high level.

In Section 2, we formally present such an estimator, suggested in Ferreira \cite{fer2}, which we shall call the runs estimator of the upcrossings. In the subsequent sections we show that such an estimator is typically weakly consistent and asymptotically normal. In Section 6 we carry out a simulation study of finite sample behaviour of the proposed estimator in a max-autoregressive process and in a first order autoregressive process. In Section 7 the performance of the estimator is also assessed through case studies in the fields of environment and finance. Conclusions are found in Section 8.

\section{The runs estimator of the upcrossings index}\setcounter{equation}{0}

\pg Before formally defining the runs estimator of the upcrossings index, we shall start by introducing some notation that will be used throughout the paper.\vspace{0.3cm}

Let  $\widehat{N}_n$ denote the point process of non-upcrossings of
$u_n,$ followed by an upcrossing,  by the first $n$ variables of
${\bf{X}},$ that is
\begin{equation}
\widehat{N}_n(B)=\sum_{i=1}^{n-3} \indi_{\{\widetilde{N}_{i,i}=0,\ X_{i+2}\leq u_n<X_{i+3}\}}\delta_{\frac{i}{n}}(B),\qquad B
\subset [0,1],\ n\geq 1.\label{p_non}
\end{equation}

%From the previous result it is clear that when condition $\widetilde{D}^{(3)}({\bf{u}})$ holds
%this point process has an important role in the characterization of
%the compound Poisson process that arises as limite in distribution
%of the upcrossings point process.\vspace{0.3cm}

For a sequence of real levels $\{u_n\}_{n\geq 1}$ lets define the random variables  $Y_i\equiv Y_i(u_n),\
i=1,\ldots,n,$ corresponding to the number of consecutive upcrossings of the level $u_n$ occurring from instant $i+2$ on, that is,
\begin{eqnarray*}
Y_i=\left\{\begin{tabular}{ll} $0\ \textrm{ if }\ \widetilde{N}_{i+2,i+2}=0$ \\
$k\ \textrm{ if }\ X_{i+2}\leq u_n<X_{i+3},X_{i+4}\leq
u_n<X_{i+5},\ldots,X_{i+2k}\leq u_n<X_{i+2k+1},
\widetilde{N}_{i+2k+2,i+2k+2}=0,$
\end{tabular}\right.
\end{eqnarray*}
with $k\geq 1.$ Furthermore, lets denote by $Z_i(u_n),$ $i=1,\ldots,n,$ the length of each of these sequences given the occurrence of a non-upcrossing followed by an upcrossing at instant $i+2,$ of level $u_n,$ which has distribution $\widetilde{\pi}_n^*(\cdot)$ given in (\ref{pi}) since
\begin{equation}
\widetilde{\pi}_n^*(k)=P(Z_i(u_n)=k)=P(Y_i=k\ |\ \widetilde{N}_{i,i}=0,\ X_{i+2}\le u_n <
X_{i+3}), \ k\geq 1,\label{dist}
\end{equation}
independent of $i$ from the stationarity of ${\bf{X}}.$

If ${\bf{X}}$ has upcrossings index $\eta>0$ and the conditions of
Proposition 1.2 hold then for $u_n=u_n^{(\nu)}$
\begin{eqnarray*}
\frac{1}{\eta}&=&\lim_{n\to +\infty}E[Z_1(u_n)]=\lim_{n\to
+\infty}E[Y_1\ |\ \widetilde{N}_{1,1}=0,X_{3}\le u_n < X_{4}]\\[0.3cm]
&=&\lim_{n\to +\infty} \frac{P(X_1\leq u_n
<X_2)}{P(\widetilde{N}_{1,1}=0,\ X_{3}\leq u_n < X_{4})}=\lim_{n\to
+\infty}\frac{E[\widetilde{N}_n([0,1])]}{E[\widehat{N}_n([0,1])]}.
\end{eqnarray*}

From this result it is natural to propose the non parametric estimator for $\eta$ given by the ratio between the total number of non-uprcrossings followed by an upcrossings  and the total number of upcrossings
\begin{equation}
\widehat{\eta}_n=\widehat{\eta}_n(u):=\frac{\widehat{N}_n([0,1])}{\widetilde{N}_n([0,1])}=\frac{\sum_{i=1}^{n-3}
\indi_{\{\widetilde{N}_{i,i}=0,\ X_{i+2}\leq u<X_{i+3}\}}}{\sum_{i=1}^{n-1}\indi_{\{X_{i}\leq u<X_{i+1}\}}},\label{estim}
\end{equation} where $u$ is a suitable threshold.

We shall call this estimator the runs estimator of the upcrossings index attending to its similarity with the runs estimator of the extremal index proposed by Nandagopalan \cite{nan} for  stationary sequences satisfying condition $D''({\bf{u}}),$ as suggested in Ferreira \cite{fer2}. In practical applications we will always use this estimator (\ref{estim}), nevertheless, the theoretical properties presented in the following sections will be proved for an estimator that is asymptotically equivalent to $\widehat{\eta}_n$ and that we define in what follows.\vspace{0.3cm}

Let  $\widehat{N}_n^{(T)}$ denote the marked point process  on $[0,1]$ defined by
\begin{equation}
\widehat{N}_n^{(T)}(B)=\sum_{i=1}^{n-3}T(Y_i)\indi_{\{\widetilde{N}_{i,i}=0,\ X_{i+2}\leq u_n<
X_{i+3}\}}\delta_{\frac{i}{n}}(B),\qquad
B\subset [0,1],\ n\geq 1,\label{marked}
\end{equation}
where $T:\N\longrightarrow \R$ is a
given mapping. Thus, $\widehat{N}_n^{(T)}$ has mass equal to
$T(Y_i)$ at the point $i/n$ whenever ${\bf{X}}$ has a non-upcrossing at $i$
followed by an upcrossing.\vspace{0.3cm}

If in (\ref{marked}) we consider $T(y)\equiv 1$ we obtain the point process $\widehat{N}_n$ in (\ref{p_non}). Whereas, if we consider $T(y)=y\geq 1$ we obtain a point process, that we shall denote by $\overline{N}_n,$ which differs from the upcrossings point process $\widetilde{N}_n$ in (\ref{p_up}) on an event with probability bounded by $P(X_1\leq u_n<X_2)$ , which converges to zero, as $n\rain {},$ since the levels  $u_n$ are commonly chosen to satisfy (\ref{eq1}). Therefore, we can instead consider the estimator $\widehat{\eta}^*_n=\widehat{N}_n([0,1]) /\overline{N}_n([0,1])$ since, for such levels, $\widehat{\eta}_n=\widehat{\eta}^*_n+o_p(1).$  The properties proved, in what follows, for $\widehat{\eta}^*_n,$ then also apply to $\widehat{\eta}_n$ that we shall use in the simulation studies.

\section{Weak consistency} \setcounter{equation}{0}

\pg We show, in this section, that $\widehat{\eta}_n^*$ is a weakly
consistent estimator of $\eta$ under mild assumptions.

Throughout it will be assumed that ${\bf{X}}$ is a stationary sequence satisfying condition $\widetilde{D}^{(3)}({\bf{u}})$ and has upcrossings index $\eta>0.$ \vspace{0.3cm}

Note that when the level $u_n$ satisfies (\ref{eq1}) there are insufficient upcrossings to give statistical ``consistency'' for the estimator $\widehat{\eta}^*_n.$ That is, as $n$ increases the value of $\widehat{\eta}^*_n$ does not necessarily converge appropriately to the value $\eta.$ Nevertheless, consistency can be achieved by the use of somewhat lower levels. We shall therefore consider non-normalized levels, in the sense of (\ref{eq1}),  $v_n=u_{[n/c_n]}^{(\nu)},$ for some fixed $\nu>0,$ that satisfy
\begin{equation}
nP(X_1\leq v_n< X_2)-c_n\nu\xrightarrow
[n\rain]{}0\label{nivel}
\end{equation}
where $\{c_n\}_{n\geq 1}$ and $\{k_n\}_{n\geq 1}$ are sequences of
real levels such that $c_n,\ k_n\xrightarrow [n\rain]{}+\infty$ and
$c_n/k_n\xrightarrow [n\rain]{}0.$ Note that for this sequence of
levels we also have $\lim_{n\to
+\infty}E[Z_1(v_n)]=1/\eta.$\vspace{0.3cm}

With the same arguments used in Nandagopalan \cite{nan} we obtain the following lemma which is essential to obtain the properties of the estimator $\widehat{\eta}_n^*$ that we further present.

\begin{lema}
Suppose $\{k_n\}_{n\geq 1}$ is  a sequence of integers such that $k_n\xrightarrow [n\rain]{}+\infty$ and that there exists a sequence $\{l_n\}_{n\geq 1}$ for which
\begin{equation}
k_n[\alpha_{n,l_n-2} +P(\widetilde{N}_n([0,l_n/n])>0)]\xrightarrow [n\rain]{} 0.\label{cond}
\end{equation}
Then $$E[e^{ia_n\widehat{N}_n^{(T)}}(J_n)]-\prod_{j=1}^{k_n}E[e^{ia_n\widehat{N}_n^{(T)}}(J_{nj})]\xrightarrow [n\rain]{} 0,$$ for any sequence of real numbers $\{a_n\}_{n\geq 1},$ where $J_n\subset [0,1],\ n\geq 1,$ is a sequence of intervals such that for each $n,$ $J_n\supset \bigcup_{j=1}^{k_n} J_{nj},$ with $J_{nj},$ $j=1,\ldots,k_n$ disjoint subintervals satisfying $m(J_n)-m(\bigcup_{j=1}^{k_n} J_{nj})\leq k_n/n$ ($m(\cdot)$ denoting Lebesgue measure).
\end{lema}

\begin{nota}
Condition (\ref{cond}) is satisfied for any sequence of normalized levels $\bf{u}^{(\nu)},$ as in (\ref{eq1}), if condition $\Delta(\bf{u}^{(\nu)})$ holds, since this condition implies that $k_n\alpha_{n,l_n-2}\xrightarrow [n\rain]{} 0$ and, for levels $\bf{u}^{(\nu)},$  $k_nP(\widetilde{N}_n([0,l_n/n])>0)\xrightarrow [n\rain]{} 0.$
\end{nota}

Henceforth, we shall write $\overline{N}_n=\overline{N}_n([0,1])=\sum_{i=1}^{n-3}Y_i\indi_{\{N_{i,i}=0,\ X_{i+2}\leq u_n<
X_{i+3}\}}$ and $\widehat{N}_n=\widehat{N}_n([0,1])=\sum_{i=1}^{n-3}\indi_{\{N_{i,i}=0,\ X_{i+2}\leq u_n<
X_{i+3}\}}$ so that $\overline{N}_n$ and $\widehat{N}_n$ will now denote random variables rather than point processes, and accordingly $\overline{N}_{r_n}=\overline{N}_n([0,r_n/n])$ and $\widehat{N}_{r_n}=\widehat{N}_n([0,r_n/n]).$

Lets consider the compound Poisson random variable $$\overline{N}_n^*=\sum_{j=1}^{\widehat{N}_n^*}Z_j^*(v_n)$$ where $\widehat{N}_n^*$ denotes a Poisson random variable with mean $c_n\eta \nu$ and the random variables $Z_j^*(v_n)$ are independent and identically distributed with the same distribution as $Z_1(v_n),$ that is, distribution $\pi_n^*$ given in (\ref{dist}).

In the following results we prove, with similar arguments to the ones used in Nandagopalan \cite{nan}, that under suitable conditions the limiting distributions of $\overline{N}_n$ and $\widehat{N}_n$ are identical to, respectively, those of $\overline{N}_n^*$ and $\widehat{N}_n^*.$

\begin{teo}
Let $\{v_n\}_{n\geq 1}$ be a sequence of levels satisfying (\ref{nivel}) and suppose there exists a sequence $\{l_n\}_{n\geq 1}$ for which (\ref{cond}) holds for such levels.

If
\begin{eqnarray}
E\left[Z_j(v_n)
\indi_{\{\widehat{N}_{r_n}>1\}}\right]=o(1)\label{hip1}
\end{eqnarray}
 uniformly in $j=1,\ldots,r_n,$ ($r_n=[n/k_n]$), then
\begin{eqnarray}
E\left[\exp(itc_n^{-1}\overline{N}_n)\right]-E\left[\exp(itc_n^{-1}\overline{N}^*_n)\right]
\xrightarrow [n\rain]{}0\label{teo11}
\end{eqnarray}\vspace{-0.2cm}
and
\begin{eqnarray}
E\left[\exp(itc_n^{-1}\widehat{N}_n)\right]-E\left[\exp(itc_n^{-1}\widehat{N}^*_n)\right]
\xrightarrow [n\rain]{}0,\label{teo12}
\end{eqnarray}
for each $t\in \R.$
\end{teo}
\bdem
 To prove (\ref{teo11}), let $\Psi_{Z_1}(t)=E[\exp(itZ_1)]$ denote the characteristic function of the random variable $Z_1$ with distribution $\widetilde{\pi}_n^*(\cdot)$ given in (\ref{dist}). Since $\overline{N}^*_n$ is a compound Poisson variable we have $$E\left[\exp(itc_n^{-1}\overline{N}^*_n)\right]=\exp(-c_n\nu\eta(1-\Psi_{Z_1}(tc_n^{-1})))$$ and hence it suffices to show that
 \begin{eqnarray}
 E\left[\exp(itc_n^{-1}\overline{N}_n)\right]=\exp(-c_n\nu\eta(1-\Psi_{Z_1}(tc_n^{-1})))+o(1). \label{A}
 \end{eqnarray}
 Nevertheless, from Lemma 3.1, we can write $$E\left[\exp(itc_n^{-1}\overline{N}_n)\right]=(E\left[\exp(itc_n^{-1}\overline{N}_{r_n})\right])^{k_n}+o(1).$$

 Now following the steps of the proof of Proposition 5.3.1 of Nandagopalan \cite{nan} we obtain $$E\left[\exp(itc_n^{-1}\overline{N}_{r_n})\right]=1-\frac{c_n}{k_n}(\eta\nu+o(1))\times(1-\Psi_{Z_1}(tc_n^{-1}))+o(k_n^{-1}),$$ since (\ref{hip1}) holds and $\frac{r_n}{c_n}k_nP(N_{1,1}=0,\ X_3\leq v_n< X_4)=\eta\nu +o(1)$ as a consequence of condition $\widetilde{D}^{(3)}(v_n).$ Moreover, since $|1-\Psi_{Z_1}(tc_n^{-1})|\leq |tc_n^{-1}|E[Z_1(v_n)]$ and $E[Z_1(v_n)]\xrightarrow [n\rain]{}\eta^{-1}$ we have $$E\left[\exp(itc_n^{-1}\overline{N}_n)\right]=\left(1-\frac{c_n}{k_n}\eta\nu(1-\Psi_{Z_1}(tc_n^{-1}))+o(k_n^{-1})\right)^{k_n}+o(1),$$ from which results (\ref{A}) and thus proves (\ref{teo11}).

 Convergence (\ref{teo12}) can be established with similar arguments. \edem

\begin{nota}
If in the previous result we consider normalized levels $u_n^{(\nu)}$ then both convergence (\ref{teo11}) and (\ref{teo12}) would hold without the need of the normalizing constant $c_n,$ being, in this case, $\widehat{N}_n^*$ a  Poisson random variable with mean $\eta \nu.$
\end{nota}

The next result justifies the need to consider lower levels $v_n$ satisfying (\ref{nivel}) in order to guarantee the consistency of the estimator.

\begin{teo}
Suppose $\{v_n\}_{n\geq 1}$ is a sequence of levels satisfying (\ref{nivel}),
\begin{eqnarray}
E\left[Z_1^*(v_n) \indi_{\{Z_1^*(v_n)>c_n\}}\right]
\xrightarrow [n\rain]{}0\label{hip2}
\end{eqnarray}
and
\begin{eqnarray}
E\left[(Z_1^*(v_n))^2
\indi_{\{Z_1^*(v_n)\leq c_n\}}\right] \xrightarrow [n\rain]{}0,\label{hip3}
\end{eqnarray} then
$$c_n^{-1}\overline{N}_n^* \xrightarrow [n\rain]{P} \nu.$$
\end{teo}
\bdem Start by noting that convergence (\ref{hip2}) along with the fact that $E[Z_1(v_n)]\xrightarrow [n\rain]{} \eta^{-1}$ implies that $\eta \nu E[Z_1^*(v_n)\indi_{\{Z_1^*(v_n)\leq c_n\}}]\xrightarrow [n\rain]{} \nu.$ It is then sufficient to show that
$$c_n^{-1}\overline{N}_n^*-\eta \nu E[Z_1^*(v_n)\indi_{\{Z_1^*(v_n)\leq c_n\}}]\xrightarrow [n\rain]{P} 0.$$

Now evoking similar arguments to the ones used by Nandagopalan \cite{nan} (Proposition 5.3.2) and considering $\gamma_n=E[Z_1^*(v_n)\indi_{\{Z_1^*(v_n)\leq c_n\}}]$, we obtain, for any $\epsilon>0$ and $\delta <\epsilon/\max\{\gamma_n\},$
\begin{eqnarray*}
\lefteqn{P(|c_n^{-1}\overline{N}_n^*-\eta\nu \gamma_n|>2\epsilon)}\\[0.3cm]
&\leq& (\delta+\eta\nu)E[Z_1^*(v_n)\indi_{\{Z_1^*(v_n)> c_n\}}] +\frac{E[(Z_1^*(v_n)\indi_{\{Z_1^*(v_n)\leq c_n\}})^2]}{c_n(\eta\nu-\delta)+(\epsilon/(\delta+\eta\nu))^2}+P(|c_n^{-1}\widehat{N}_n^*-c_n\eta\nu|>\delta c_n).
\end{eqnarray*}

The first two terms tend to zero by (\ref{hip2}) and (\ref{hip3}) respectively, while the third term tends to zero since $\widehat{N}_n^*$ is a compound Poisson random variable with mean $c_n\eta \nu.$ Hence, since $\epsilon>0$ is arbitrary, we obtain $c_n^{-1}\overline{N}_n^*-\eta\nu \gamma_n\xrightarrow [n\rain]{P} 0$ as required.\edem

We proved in Theorem 3.1 that $c_n^{-1}\overline{N}_n$ and $c_n^{-1}\widehat{N}_n$ have the same asymptotic distribution as, respectively,  $c_n^{-1}\overline{N}_n^*$ and $c_n^{-1}\widehat{N}_n^*,$  and in Theorem 3.2 proved that $c_n^{-1}\overline{N}_n^*\xrightarrow [n\rain]{P} \nu.$ Moreover, since $\widehat{N}_n^*$ is a Poisson random variable with mean $c_n\eta \nu$ it follows that $c_n^{-1}\widehat{N}_n^*\xrightarrow [n\rain]{P} \eta\nu.$ The weak consistency of $\widehat{\eta}^*_n$ is now an immediate consequence, summarized in the next result.

\begin{coro}
If the conditions of Theorems 3.1 and 3.2 hold then  $$\widehat{\eta}_n^* \xrightarrow [n\rain]{P} \eta.$$
\end{coro}

\section{Asymptotic Normality} \setcounter{equation}{0}

\pg Imposing additional  conditions on the limiting behaviour of the first and second moments of the variable  $Z_1(v_n)$ we obtain in this section the asymptotic normality of our estimator $\widehat{\eta}_n^*.$

The proof of the next result, which is essential in obtaining what follows, shall be omitted since it is similar to the proof of Proposition 5.4.1 of Nandagopalan \cite{nan}.

\begin{teo}
Suppose that $\{v_n\}_{n\geq 1}$ is  a sequence of levels satisfying (\ref{nivel}) and suppose there exists a sequence $\{l_n\}_{n\geq 1}$ for which (\ref{cond}) holds for such levels.

If
$$E[\overline{N}_{r_n}^2\indi_{\{\widehat{N}_{r_n}>1\}}]=o(c_n/k_n),
$$
\begin{eqnarray}
E[(Z_j(v_n))^2\indi_{\{\widehat{N}_{r_n}>1\}}]=o(1)\label{hip4}
\end{eqnarray} uniformly in $j=1,\dots,
r_n,$ $(r_n=[n/k_n])$ and $$\sigma^2_n=E[(Z_1(v_n))^2]$$ is a bounded sequence, then $\displaystyle c_n^{-1}\left[\begin{tabular}{l}
$\overline{N}_n-E[\overline{N}_n]$
\\ $\widehat{N}_n-E[\widehat{N}_n]$\end{tabular}\right]$
converges in distribution if and only if  $\displaystyle
c_n^{-1}\left[\begin{tabular}{l} $\overline{N}^*_n-E[\overline{N}^*_n]$ \\
$\widehat{N}^*_n-E[\widehat{N}^*_n]$\end{tabular}\right]$ does and, in this case, the limits coincide.
\end{teo}

\begin{teo}\label{Art153prop542}
Suppose that $\{v_n\}_{n\geq 1}$ is  a sequence of levels satisfying (\ref{nivel}) and suppose there exists a sequence $\{l_n\}_{n\geq 1}$ for which (\ref{cond}) holds for such levels.

If
\begin{eqnarray}
\sigma^2_n=E[(Z_1(v_n))^2]\xrightarrow [n\rain]{} \sigma^2<+\infty\label{hip5}
\end{eqnarray}
and for each $\epsilon>0,$
\begin{eqnarray}
E[(Z_1^*(v_n))^2\indi_{\{(Z_1^*(v_n))^2>\epsilon c_n\}}]\xrightarrow [n\rain]{} 0.\label{hip6}
\end{eqnarray}
Then $$c_n^{-1/2}\left[\begin{tabular}{l} $\overline{N}_n^*-E[\overline{N}^*_n]$ \\
$\widehat{N}^*_n-E[\widehat{N}^*_n]$
\end{tabular}\right]\xrightarrow [n\rain]{d} \mathcal{N}
\left[\left[\begin{tabular}{l} $0$\\
$0$\end{tabular}\right],\left[\begin{tabular}{cc} $\nu\eta\sigma^2$ & $\nu$
\\ $\nu$ & $\nu\eta$ \end{tabular}\right]\right].
$$
\end{teo}
\bdem  Since $\overline{N}_n^*$ is a compound Poisson random variable it holds $E[\overline{N}_n^*]=E[\widehat{N}_n^*]E[Z_1(v_n)]$ and we can write
\begin{equation}
\overline{N}^*_n-E[\overline{N}^*_n]=(\overline{N}^*_n-\widehat{N}_n^*E[Z_1(v_n)])+E[Z_1(v_n)](\widehat{N}_n^*-E[\widehat{N}^*_n]).\label{Art153545}
\end{equation}
From Lindeberg's condition (\ref{hip6}) the central limit theorem holds, this is,  $$
 c_n^{-1/2}(\overline{N}^*_n-\widehat{N}_n^*E[Z_1(v_n)])\xrightarrow [n\rain]{d} \mathcal{N}\left(0,
 \eta\nu\left(\sigma^2-\frac{1}{\eta^2}\right)\right)$$
and therefore
 \begin{eqnarray*}
 E[\exp(itc_n^{-1/2}(\overline{N}^*_n-\widehat{N}_n^*E[Z_1(v_n)]))]&=&\exp(-c_n\eta\nu[1-\Psi_{Z_1}(tc_n^{-1/2})e^{-itc_n^{-1/2}E[Z_1(v_n)]}])
 \nonumber\\[0.3cm]
    &\xrightarrow [n\rain]{} & \exp\left(\frac{\eta\nu}{2}\left(\sigma^2-\frac{1}{\eta^2}\right)t^2\right).
\end{eqnarray*}
%and consequently
%\begin{eqnarray}
%-c_n\eta\nu[1-\Psi_{Z_1}(tc_n^{-1/2})e^{-itc_n^{-1/2}E[Z_1(v_n)]}]\xrightarrow [n\rain]{}\frac{\eta\nu}{2}\left(\sigma^2-\frac{1}{\eta^2}\right)t^2.   \label{eq1}
%\end{eqnarray}

Also,  $c_n^{-1/2}(\widehat{N}_n^*-E[\widehat{N}^*_n])\xrightarrow [n\rain]{d} \mathcal{N}(0, \eta\nu),$
so that $$
  E[\exp\{isc_n^{-1/2}(\widehat{N}_n^*-E[\widehat{N}^*_n])\}]=\exp(-c_n\eta\nu(1-e^{isc_n^{-1/2}}+isc_n^{-1/2}))
  \xrightarrow [n\rain]{}\exp\left(-\frac{\eta\nu s^2}{2}\right).$$
Now \begin{eqnarray}
 \lefteqn{E[\exp\{itc_n^{-1/2}(\overline{N}^*_n-\widehat{N}_n^*E[Z_1(v_n)])+isc_n^{-1/2}(\widehat{N}_n^*-E[\widehat{N}^*_n])\}]}\nonumber\\[0.3cm]
 &=&\exp(-c_n\eta\nu(1-\Psi_{Z_1}(tc_n^{-1/2})e^{-itc_n^{-1/2}E[Z_1(v_n)]}e^{isc_n^{-1/2}}+isc_n^{-1/2})),\label{Art153548}
\end{eqnarray}
with $\Psi_{Z_1}(tc_n^{-1/2})=1+itc_n^{-1/2}E[Z_1(v_n)]+\rho_n^{(1)}$ and $e^{-itc_n^{-1/2}E[Z_1(v_n)]}=1-itc_n^{-1/2}E[Z_1(v_n)]+\rho_n^{(2)}$ where $\rho_n^{(1)}$ and $\rho_n^{(2)}$ are bounded by $t^2c_n^{-1}\sigma^2_n/2.$

Thus,
\begin{eqnarray*}
\lefteqn{\eta\nu c_n[1-\Psi_{Z_1}(tc_n^{-1/2})e^{-itc_n^{-1/2}E[Z_1(v_n)]}e^{isc_n^{-1/2}}+isc_n^{-1/2}]}\\[0.3cm]
&=&\eta\nu c_n[(1-\Psi_{Z_1}(tc_n^{-1/2})e^{-itc_n^{-1/2}E[Z_1(v_n)]})+(1-e^{isc_n^{-1/2}}+isc_n^{-1/2})+o(c_n^{-1})]\\[0.3cm]
&\xrightarrow [n\rain]{}&\frac{\eta\nu(\sigma^2-1/\eta^2)t^2}{2}+\frac{\eta\nu s^2}{2}.
\end{eqnarray*}

This along with (\ref{Art153548}) implies that
\begin{equation}
c_n^{-1/2}\left[\begin{tabular}{l} $\overline{N}_n^*-\widehat{N}_n^*E[Z_1(v_n)]$ \\
$\widehat{N}^*_n-E[\widehat{N}^*_n]$ \end{tabular}\right]\xrightarrow [n\rain]{d} \mathcal{N} \left(\left[\begin{tabular}{l} $0$\\
$0$\end{tabular}\right],\left[\begin{tabular}{cc} $\eta\nu\left(\sigma^2-\frac{1}{\eta^2}\right)$ & $0$ \\
$0$ & $\nu\eta$
\end{tabular}\right]\right),\nonumber
\end{equation}
which together with (\ref{Art153545}) and the fact that  $E[Z_1(v_n)]\xrightarrow [n\rain]{} 1/\eta$  yields the result. \edem

The asymptotic normality of $\widehat{\eta}^*_n$ is now an immediate consequence of the previous results.

\begin{coro}
 If the conditions of Theorems 4.1 and 4.2 hold then
 $$\sqrt{c_n}(\widehat{\eta}^*_n-\eta_n)\xrightarrow [n\rain]{d}
 \mathcal{N}\left(0,\sqrt{\frac{\eta}{\nu}(\eta^2\sigma^2-1)}\right),$$
 where $\eta_n=E[\widehat{N}_n]/E[\overline{N}_n].$
\end{coro}
\bdem Conditions (\ref{hip4})-(\ref{hip6}) imply the conditions of Theorems 3.1 and 3.2, hence $c_n^{-1}\overline{N}_n\xrightarrow [n\rain]{P}\nu.$ On the other hand condition $\widetilde{D}^{(3)}(v_n)$  implies that $\eta_n\xrightarrow [n\rain]{}\eta.$ The result now follows from  the fact that
$$
 \sqrt{c_n}(\widehat{\eta}_n-\eta_n)=\frac{1}{c_n^{-1}\overline{N}_n}\{c_n^{-1/2}(\widehat{N}_n-E[\widehat{N}_n])-\eta_nc_n^{-1/2}(\overline{N}_n-
 E[\overline{N}_n])\}.$$\edem

\begin{nota}
Since the variance of $\widehat{\eta}^*_n$ is of order $1/\nu,$ if we kept $\nu$ fixed as in (\ref{eq1}) we could not guarantee the consistency of $\widehat{\eta}^*_n.$ Hence the need to assume that $\nu=\nu_n\xrightarrow []{} +\infty$ as $n \xrightarrow []{} +\infty.$
\end{nota}

 To conclude we present a result that enables one to construct approximate  confidence intervals or a hypothesis test regarding $\eta,$ {\it{i.e.}}, to determine the extent of clustering of upcrossings of high levels in the observed data. From a practical viewpoint it is more useful than  Corollary 4.1.

 \begin{coro}
 Suppose that the conditions of Corollary 4.1 hold, $$\sqrt{c_n}(\eta_n-\eta)\xrightarrow [n\rain]{}0$$ and $$c_n^{-1}E[(Z_1^*(v_n))^4\indi_{\{(Z_1^*(v_n))^2\leq c_n\}}]\xrightarrow [n\rain]{}0.$$
 Then $$\sqrt{\frac{\overline{N}_n}{\widehat{\eta}_n^*((\widehat{\eta}_n^{*}\widehat{\sigma}_n)^2-1)}}\ (\widehat{\eta}_n^*-\eta)\xrightarrow [n\rain]{} \mathcal{N} (0,1),$$ where $\widehat{\sigma}_n^2=\frac{\sum_{i=1}^{n-3}Y_i^2\indi_{\{\widetilde{N}_{i,i}=0,\ X_{i+2}\leq v_n<X_{i+3}\}}}{\widehat{N}_n}.$
 \end{coro}
 \bdem Straightforward from Corollary 4.1 and the fact that Theorems 3.1 and 3.2 imply that $$c_n^{-1}\sum_{i=1}^{n-3}Y_i^2\indi_{\{\widetilde{N}_{i,i}=0,\ X_{i+2}\leq v_n<X_{i+3}\}}\xrightarrow [n\rain]{P} \eta \nu \sigma^2.$$
 \edem

We recall that the properties proved for the estimator $\widehat{\eta}_n^*$ remain valid for the estimator $\widehat{\eta}_n$ used in the sequel of this paper.

\section{The choice of the levels}

\pg In Section 3 it was shown that for the runs estimator of the upcrossings index $\widehat{\eta}_n$ consistency can be achieved by the use of somewhat lower levels $v_n$ satisfying (\ref{nivel}). Thus, the precise choice of $v_n$ depends on the knowledge of the joint distribution of $(X_1,X_2),$ typically unknown. These deterministic levels will have to be replaced, in practical situations, by random levels suggested by the relation
\begin{equation}
\frac{n}{c_n}P(X_1\leq v_n<X_2)\sim \nu.\label{norm}
 \end{equation}
This relation basically states  that the expected number of upcrossings of the level $v_n$ is approximately $c_n \nu.$  Contrarily to the random levels used in the estimation of the extremal index, these random levels can not be represented by an appropriate order statistic. Nevertheless, for levels $v_n$ such that $\frac{n}{c_n}P(X_1>v_n)\sim \tau$ the expected number of upcrossings will be necessarily smaller or equal to $c_n\tau.$
It then seems natural to also replace $v_n$ by the appropriate order statistic used in the estimation of the extremal index (for example, Nandagopalan \cite{nan}), namely $\widehat{v}_n=X_{n-[c_n\tau]:n}.$ With this consideration we obtain the following estimator for $\eta$ $$\widetilde{\eta}_n=\frac{\widehat{N}_n(\widehat{v}_n)}{\widetilde{N}_n(\widehat{v}_n)},$$ where $\widehat{N}_n(\widehat{v}_n)=\widehat{N}_n([0,1])$ and  $\widetilde{N}_n(\widehat{v}_n)=\widetilde{N}_n([0,1])$ with $u_n=\widehat{v}_n.$

The weak consistency of this estimator can be obtained by showing that it is closely approximated by the corresponding estimates based on the non-random levels $v_n$. For this, lets start by noting that for two levels $v_1$ and $v_2,$ such that $v_1\neq v_2,$ and fixed $n$ we have
\begin{equation}
|\widetilde{N}_n(v_1)-\widetilde{N}_n(v_2)|\leq |N_n(v_1)-N_n(v_2)|\label{niv1}
\end{equation}
and
\begin{equation}
|\widehat{N}_n(v_1)-\widehat{N}_n(v_2)|\leq |N_n(v_1)-N_n(v_2)|\label{niv2}
\end{equation}
where $N_n(v_n)=\sum_{i=1}^n \indi_{\{X_i>v_n\}}.$

\begin{teo}
Suppose that for each $\nu>0$ there exists $u_{[n/c_n]}^{(\nu)}=u_{[n/c_n]}^{(\tau)}$ for some $\tau>0,$ the conditions of Proposition 5.3.1 and 5.3.2 in Nandagopalan \cite{nan} hold for each $\tau'$ in a neighbourhood of $\tau$ and conditions (3.3), (3.7) and (3.8) hold for each $\nu'$ in a neighbourhood of $\nu.$ Then $\widetilde{\eta}_n \xrightarrow [n\rain]{P} \eta.$
\end{teo}

\bdem
For $\epsilon>0$
\begin{eqnarray}
\lefteqn{P(|c_n^{-1}\widetilde{N}_n(\widehat{v}_n)-\nu|>6\epsilon)}\nonumber\\
&\leq& P(|c_n^{-1}[c_n\tau]-N_n(v_n)|>3\epsilon)+P(|c_n^{-1}\widetilde{N}_n(v_n)-\nu|>3\epsilon)\label{conv1}
\end{eqnarray}
by (\ref{niv1}), with $v_1$ and $v_2$ replaced by $\widehat{v}_n$ and $v_n=u_{[n/c_n]}(\tau+\epsilon).$ If $\epsilon$ is sufficiently small, then it follows from the results in \cite{nan} and in Section 3 that $c_n^{-1}N_n(v_n)\xrightarrow [n\rain]{P} \tau+\epsilon$ and $c_n^{-1}\widetilde{N}_n(v_n)\xrightarrow [n\rain]{P} \nu+\epsilon.$ Therefore, since $c_n^{-1}[c_n\tau]\xrightarrow [n\rain]{}\tau,$ for large $n$ (\ref{conv1}) is dominated by $$P(|c_n^{-1}N_n(v_n)-(\tau+\epsilon)|>\epsilon)+P(||c_n^{-1}\widetilde{N}_n(v_n)-(\nu+\epsilon)|>\epsilon)\xrightarrow [n\rain]{}0,$$ which proves that $c_n^{-1}\widetilde{N}_n(\widehat{v}_n)\xrightarrow [n\rain]{P}\nu.$
On the other hand, for $\epsilon>0,$ we have from (\ref{niv2}), with $v_1$ and $v_2$ replaced by $\widehat{v}_n$ and $v_n=u_{[n/c_n]}(\tau+\epsilon),$ that
\begin{eqnarray}
\lefteqn{P(|c_n^{-1}\widehat{N}_n(\widehat{v}_n)-\eta\nu|>6\epsilon)}\nonumber\\
&\leq& P(|c_n^{-1}[c_n\tau]-N_n(v_n)|>3\epsilon)+P(|c_n^{-1}\widehat{N}_n(v_n)-\nu|>3\epsilon)\label{conv2}.
\end{eqnarray}
Hence, if $\epsilon$ is sufficiently small we have  $c_n^{-1}\widehat{N}_n(v_n)\xrightarrow [n\rain]{P} \eta(\nu+\epsilon)$ (results of Section 3) and can may conclude that (\ref{conv2}) is dominated by $$P(|c_n^{-1}N_n(v_n)-(\tau+\epsilon)|>\epsilon)+P(||c_n^{-1}\widehat{N}_n(v_n)-\eta(\nu+\epsilon)|>\epsilon)\xrightarrow [n\rain]{}0,$$ which proves that $c_n^{-1}\widehat{N}_n(\widehat{v}_n)\xrightarrow [n\rain]{P}\eta\nu.$

The result is now an immediate consequence of the two convergences established.\edem

The proof of the asymptotic normality of $\widetilde{\eta}_n$ remains an open problem.

\section{Simulations: Some examples} \setcounter{equation}{0}

\pg This section studies some examples of sequences exhibiting clustering of upcrossings. We illustrate the performance of the runs estimator of the upcrossings index with the max-autoregressive (ARMAX) process of Ferreira \cite{fer1} and the negatively correlated first order autoregressive (AR(1)) process of Chernick {\it{et al.}} \cite{cher}, for which the values of the upcrossings index $\eta$ are well known.

\subsection{ARMAX process}

\pg Lets consider the ARMAX process of Ferreira \cite{fer1},
\begin{eqnarray}
X_n=\max\{Y_n,\ Y_{n-2},\ Y_{n-3}\}, \quad n\geq 1,\label{armax}
\end{eqnarray}
where$\{Y_n\}_{n\geq -2}$ is a a sequence of independent and uniformly distributed on $[0,1]$ variables. In Figure 1 we present a sample path of this process.

\begin{figure}[!ht]
\begin{center}
\includegraphics[scale=0.4]{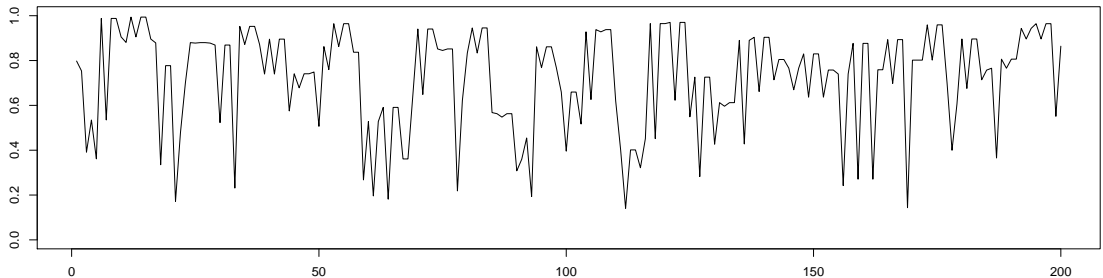}\vspace{-0.5cm}
\caption{Sample path of the stationary  ARMAX process in (\ref{armax}).}
\end{center}
\end{figure}

Condition $\widetilde{D}^{(3)}({\bf{u}})$ holds for this stationary sequence, with ${\bf{u}}={\bf{u}}^{(\nu)}=\{u_n=1-\tau/n\}_{n\geq 1},$ $\tau>0.$ It has extremal index $\theta=1/3$ and upcrossings index $\eta=1/2.$ Furthermore, it holds \linebreak $\lim_{n\to +\infty}nP(X_1>u_n)=3\tau$ and $\lim_{n\to +\infty}nP(X_1\leq u_n< X_2)=2\tau,$ so the levels $u_n$ are simultaneously normalized for exceedances and upcrossings.  Note that this implies that for a sample $(X_1,\ldots,X_n)$ of (\ref{armax}) with $n$ sufficiently large, the number of upcrossings of a high level is approximately $2/3$ of the number of exceedances of the same level.\vspace{0.3cm}

We shall con\-si\-der in the following simulations the  level $u=X_{n-k:n},$  in (\ref{estim}) corresponding to the  $(k+1)$th top order statistics associated to the random sample $(X_1,\ldots,X_n)$ of (\ref{armax}), commonly used in the estimation of the extremal index (see, for example, Nandagopalan \cite{nan} and Gomes, {\it{et al.}} \cite{gom}). Hence, the upcrossings index estimator is now a function of $k,$ $\widehat{\eta}_n(k).$

%This choice of levels and relation (\ref{rel_eta_teta}) implies that in this case we can obtain estimates for $\theta$ from the estimates of $\eta,$ %obtained from $\widehat{\eta}_n$ in  (\ref{estim}), by multiplying the latter by $2/3.$

In Figure 2 we present a sample path of the estimator $\widehat{\eta}_n$ in  (\ref{estim}) as a function of $k\geq 1,$ for a sample of size $n=5000$ in a linear scale and in a logarithmic scale.

\begin{figure}[!ht]
\begin{center}
\includegraphics[scale=0.25]{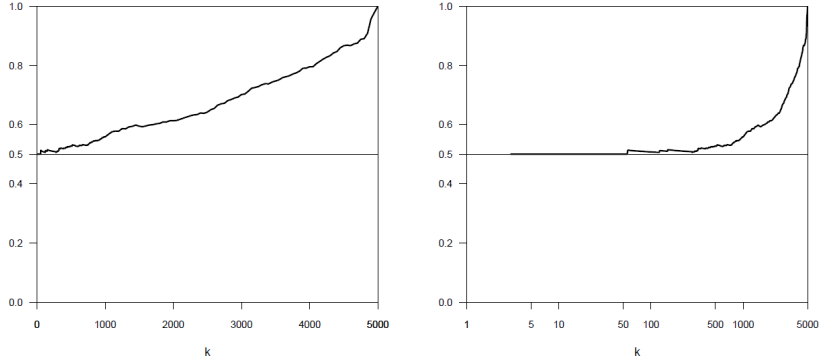}\vspace{-0.5cm}
\caption{Sample path of the upcrossings index estimator as a function of $k$ for a sample of size $n=5000$ from the stationary  ARMAX process in (\ref{armax}), in a linear scale (left) and in a logarithmic scale (right).}
\end{center}
\end{figure}

As we can see from Figure 2 the logarithmic scale enhances the performance of $\widehat{\eta}_n$ for small values of $k,$ giving better insight of its stability region around $\eta=0.5$.\vspace{0.3cm}

For samples of size $n=100, 200, 500, 1000, 2000, 5000$ and 10000, from the ARMAX process in (\ref{armax}), we have performed a multi-sample Monte Carlo simulation with 5000 runs and 10 replicates. For details on multi-sample simulation see, for instance, Gomes and Oliveira \cite{gom2}. We have simulated for the estimator $\widehat{\eta}_n(k)$ in  (\ref{estim}), the mean value ($E$), the mean squared error ($MSE$) and the optimal sample fraction $k_0$ with $k_0:=\arg \min_k MSE[\widehat{\eta}_n(k)].$

In Figure 3, we illustrate the values of the estimated mean values and MSE's of $\widehat{\eta}_n(k)$ in  (\ref{estim}), for a sample of size $n=5000$ from the ARMAX process in (\ref{armax}), with upcrossings index  $\eta=0.5.$ In Table 1  we present the main distributional properties of the estimator under study with the associated 95\% confidence intervals (see Gomes and Oliveira \cite{gom2} ).

 \begin{figure}[!ht]
\begin{center}
\includegraphics[scale=0.25]{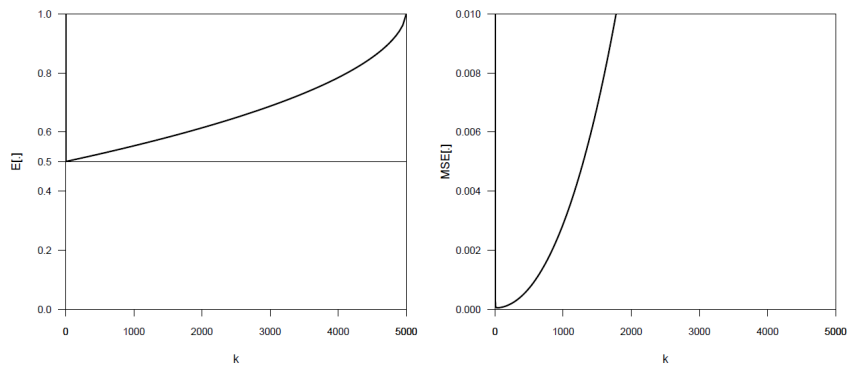}\vspace{-0.5cm}
\caption{Estimated mean values (left) and mean squared errors (right), for samples  of size $n=5000$ from the ARMAX process in (\ref{armax}) ($\eta=0.5$)}
\end{center}
\end{figure}

\begin{table}[!htb]
\centering \scriptsize
 \begin{tabular*}{\textwidth}{@{\extracolsep{\fill}}cccccc}
  \hline
 $n$ & $k_0$ & $k_0/n$ & $E[\bullet ]$ &  $MSE[\bullet ]$ & $SD[\bullet ]$ \\
  \hline
  $100$ & $16$ & $0.160$ & $0.52735 \pm 0.00042$ & $0.00497 \pm 0.000087$ & $0.06495 \pm 0.00064$ \\
  $200$ & $19$ & $0.095$ & $0.51754 \pm 0.00029$ & $0.00214 \pm 0.000029$ & $0.04282 \pm 0.00037$ \\
  $500$ & $25$ & $0.050$ & $0.50977 \pm 0.00021$ & $0.00072 \pm 0.000018$ & $0.02493 \pm 0.00030$ \\
  $1000$ & $30$ & $0.030$ & $0.50622 \pm 0.00014$ & $0.00033 \pm 0.000008$ & $0.01699 \pm 0.00022$ \\
  $2000$ & $36$ & $0.018$ & $0.50389 \pm 0.00007$ & $0.00015 \pm 0.000003$ & $0.01162 \pm 0.00013$ \\
  $5000$ & $54$ & $0.011$ & $0.50241 \pm 0.00006$ & $0.00005 \pm 0.000001$ & $0.00694 \pm 0.00010$ \\
  $10000$ & $67$ & $0.007$ & $0.50152 \pm  0.00003$ & $0.00003 \pm 0.000001$ & $0.00477 \pm 0.00006$\\
   \hline
\end{tabular*}\caption{Optimal sample fractions, mean values and mean squared errors of the estimator at its optimal levels, for the ARMAX process, with $\eta=0.5.$ }
\end{table}

\subsection{AR(1) process}

\pg Lets consider the negatively correlated uniform AR(1) process of Chernick {\it{et al.}} \cite{cher},
\begin{equation}
X_n=-\frac{1}{r}X_{n-1}+\epsilon_n,\quad n\geq 1,\label{ar}
\end{equation}
 where $\{\epsilon_n\}_{n\geq 1}$ is a sequence of independent and identically distributed random variables, such that, for a fixed integer $r\geq 2,$ $\epsilon_n\sim U\{\frac{1}{r},\ldots,\frac{r-1}{r},1\}$ and  $X_0\sim U(0,1)$ independent of $\epsilon_n.$

Condition $\widetilde{D}^{(3)}({\bf{u}})$ also holds for this stationary sequence, with ${\bf{u}}={\bf{u}}^{(\nu)}=\{u_n=1-\tau/n\}_{n\geq 1},$ $\tau>0,$ moreover $\lim_{n\to +\infty}nP(X_1>u_n)=\tau$ and $\eta=\theta=1-1/r^2$ (see Sebastião {\it{et al.}} \cite{seb}). Condition $D''({\bf{u}})$ typically doesn't hold for these sequences since they tend to oscillate rapidly near extremes.
To illustrate this characteristic  we present in Figure 4 sample paths of the negatively correlated uniform AR(1) processes for $r=2,$ $r=3$ and $r=5$ ($\eta=\theta=0.75,$  $\eta=\theta=0.89$ and $\eta=\theta=0.96,$ respectively).

\begin{figure}[!ht]
\begin{center}
\includegraphics[scale=0.4]{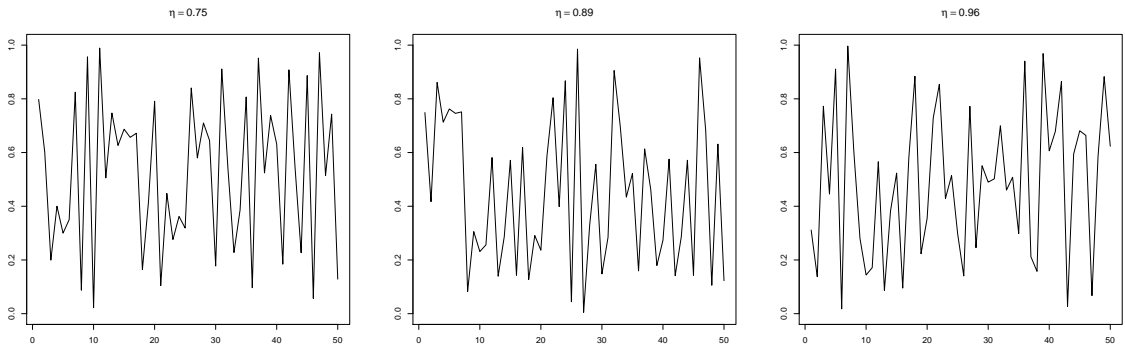}\vspace{-0.5cm}
\caption{Sample paths of the stationary  AR(1) process in (\ref{ar}), with $r=2$ (left), $r=3$ (center) and $r=5$ (right)}
\end{center}
\end{figure}

%In Figure 5 we present sample paths of the estimator $\widehat{\eta}_n$ in  (\ref{estim}) as a function of $k,$ for a sample of size 5000 of the previous process with $r=2$ and $r=3$ and $r=5$ ($\eta=0.75,\ 0.89,\ 0.96,$ respectively).

%\begin{figure}[!ht]
%\begin{center}
%\includegraphics[scale=0.4]{AR1-samplepath-n5000}\vspace{-0.5cm}
%\caption{Sample paths of the upcrossings index estimator as a function of $k$ for a sample of size $n=5000$ from the AR(1) process with $r=2$ (left), $r=3$ (center) and $r=5$ (right)}
%\end{center}
%\end{figure}

We present, in Figures 5-7, the estimated mean values and MSE's of $\widehat{\eta}_n(k)$ in  (\ref{estim}), for a sample of size $n=5000$ from the AR(1) process in (\ref{ar}), with upcrossings index  $\eta=0.75, \ 0.89$ and 0.96, respectively.

\begin{figure}[!ht]
\begin{center}
{\small{$$r=2$$}}\vspace{-1cm}

\includegraphics[scale=0.25]{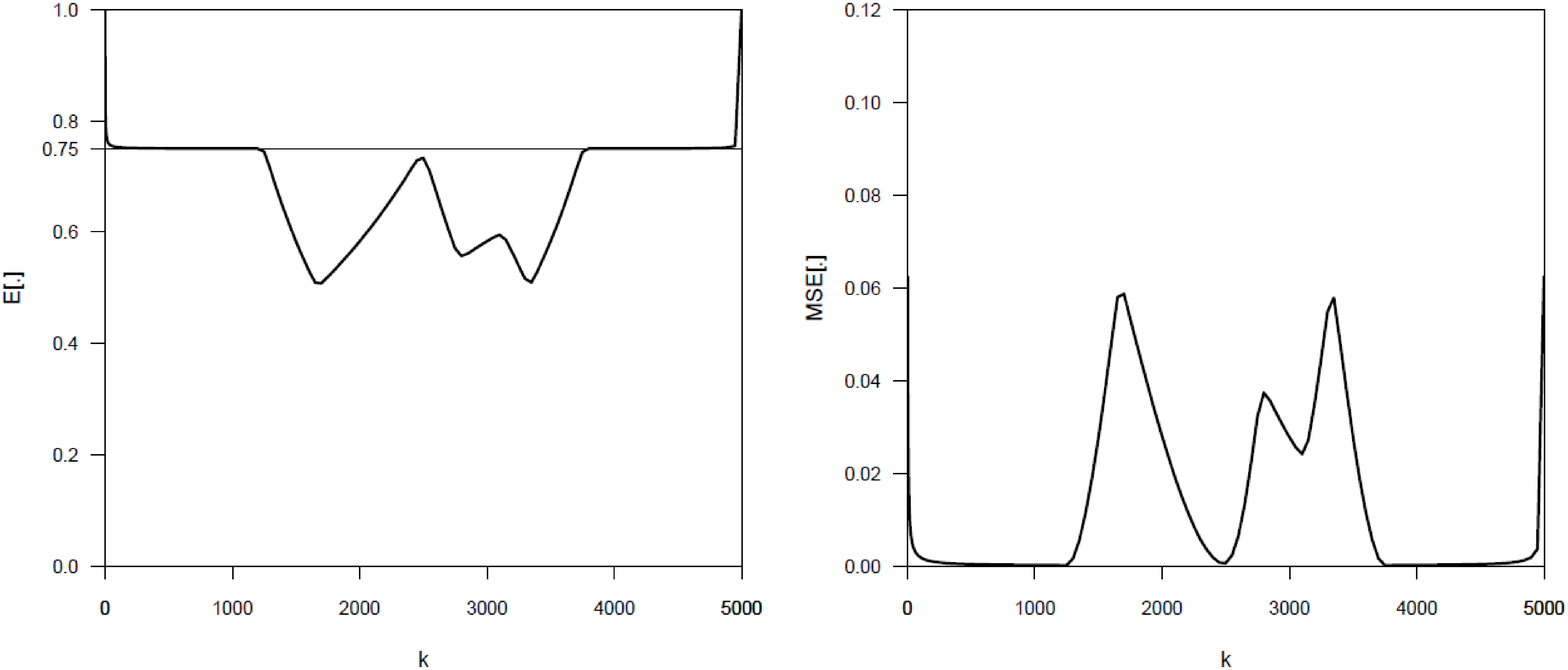}\vspace{-1cm}

$$r=3$$\vspace{-1cm}

\includegraphics[scale=0.25]{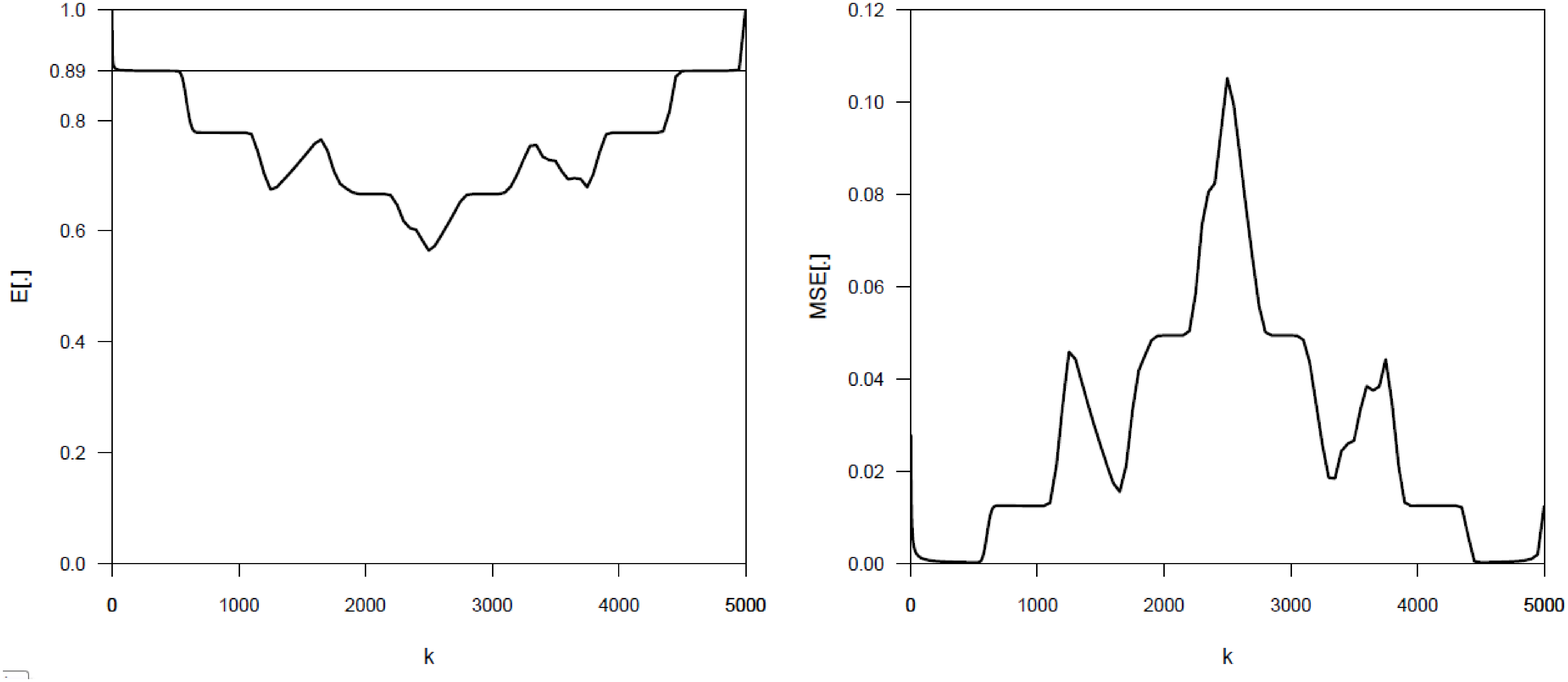}\vspace{-1cm}

$$r=5$$\vspace{-1cm}

\includegraphics[scale=0.25]{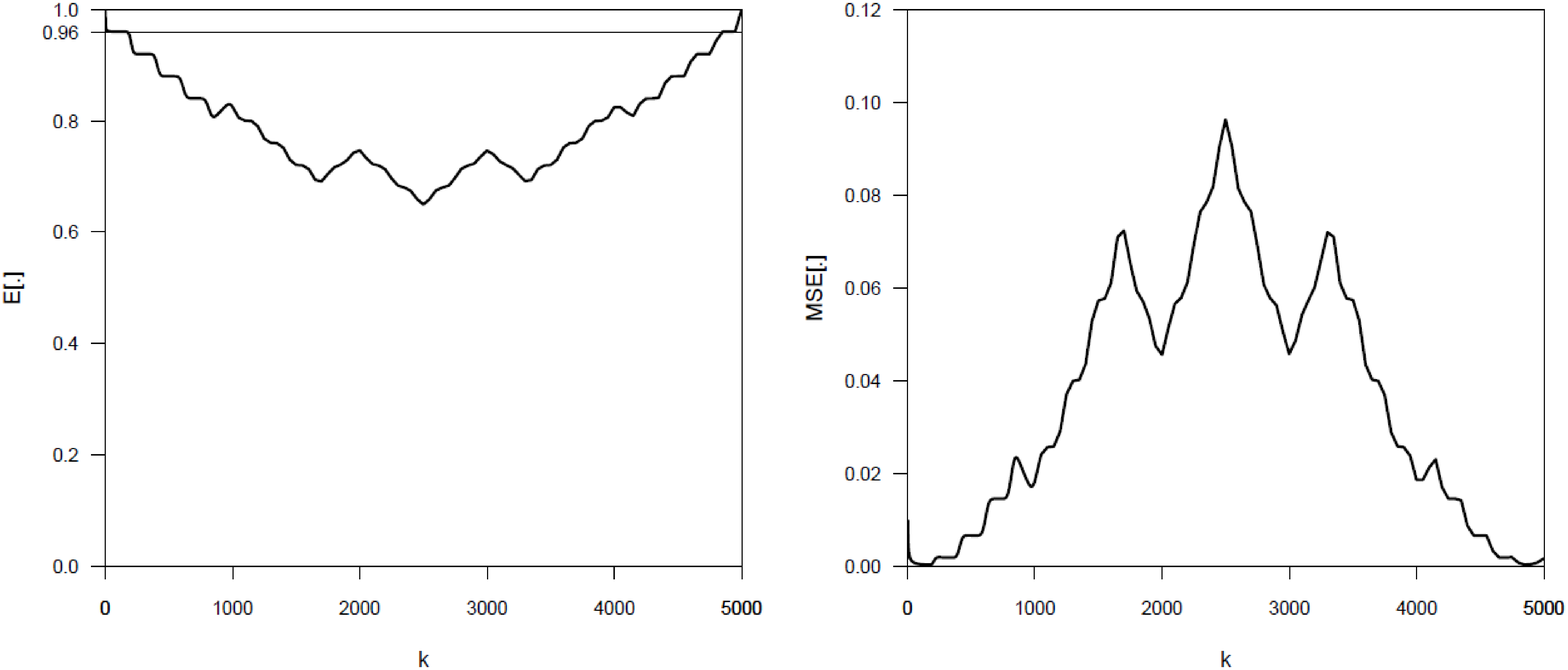}
\vspace{-0.5cm}
\caption{Estimated mean values (left) and mean squared errors (right), for samples  of size $n=5000$ from the AR(1) process in (\ref{ar}) with, respectively, $r=2$ ($\eta=0.75$), $r=3$ ($\eta=0.89$) and $r=5$ ($\eta=0.96$)}
\end{center}
\end{figure}
%\begin{figure}[!ht]
%\begin{center}
%\includegraphics[scale=0.5]{AR1-runs-n5000-m5000-r3}\vspace{-0.5cm}
%\caption{Estimated mean values (left) and mean squared errors (right), for samples  of size $n=5000$ from the AR(1) process in (\ref{ar}) with $r=3$  ($\eta=0.89$)}
%\end{center}
%\end{figure}\vspace{-0.5cm}
%\begin{figure}[!ht]
%\begin{center}
%\includegraphics[scale=0.5]{AR1-runs-n5000-m5000-r5}\vspace{-0.5cm}
%\caption{Estimated mean values (left) and mean squared errors (right), for samples  of size $n=5000$ from the AR(1) process in (\ref{ar}) with $r=5$  ($\eta=0.96$)}
%\end{center}
%\end{figure}

In Table 2 we present the main distributional properties of the estimator under study with the associated 95\% confidence intervals.

\begin{table}[!htb]
\begin{center}\scriptsize
 \begin{tabular*}{\textwidth}{@{\extracolsep{\fill}}ccccc}
  \hline
 $n$ & $k_0$ & $k_0/n$ & $E[\bullet ]$ &  $MSE[\bullet ]$  \\
  \hline
  ${\boldsymbol{\eta}}{\bf{=0.75}}$&&&&\\
  $100$  & $24$   & $0.240$ & $0.72839 \pm 0.00081$ & $0.00769 \pm 0.000075$  \\
  $200$  & $49$   & $0.245$ & $0.73228 \pm 0.00043$ & $0.00364 \pm 0.000055$  \\
  $500$  & $123$  & $0.246$ & $0.73942 \pm 0.00026$ & $0.00140 \pm 0.000016$  \\
  $1000$ & $246$  & $0.246$ & $0.74375 \pm 0.00009$ & $0.00068 \pm 0.000009$  \\
  $2000$ & $494$  & $0.247$ & $0.74576 \pm 0.00017$ & $0.00034 \pm 0.000005$  \\
  $5000$ & $1240$ & $0.248$ & $0.74745 \pm 0.00010$ & $0.00013 \pm 0.000002$  \\
  ${\boldsymbol{\eta}}{\bf{=0.89}}$&&&&\\
  $100$  & $10$  & $0.100$ & $0.86547 \pm 0.00112$ & $0.01254 \pm 0.000240$ \\
  $200$  & $19$  & $0.095$ & $0.87623 \pm 0.00077$ & $0.00598 \pm 0.000096$ \\
  $500$  & $48$  & $0.096$ & $0.88347 \pm 0.00029$ & $0.00221 \pm 0.000025$ \\
  $1000$ & $98$  & $0.098$ & $0.88642 \pm 0.00017$ & $0.00105 \pm 0.000009$ \\
  $2000$ & $201$ & $0.101$ & $0.88768 \pm 0.00020$ & $0.00050 \pm 0.000007$ \\
  $5000$ & $519$ & $0.104$ & $0.88842 \pm 0.00007$ & $0.00019 \pm 0.000003$ \\
  ${\boldsymbol{\eta}}{\bf{=0.96}}$&&&&\\
  100 & 7 & 0.070    & $0.91442 \pm 0.00123$ & $0.01431 \pm 0.000366$ \\
  200 & 9 & 0.045    & $0.93962 \pm 0.00081$ & $0.00728 \pm 0.000126$ \\
  500 & 17 & 0.034   & $0.95373 \pm 0.00072$ & $0.00274 \pm 0.000077$ \\
  1000 & 33 & 0.033  & $0.95687 \pm 0.00054$ & $0.00130 \pm 0.000038$ \\
  2000 & 67 & 0.034  & $0.95842 \pm 0.00026$ & $0.00061 \pm 0.000011$ \\
  5000 & 173 & 0.035 & $0.95949 \pm 0.00017$ & $0.00023 \pm 0.000004$ \\
   \hline
\end{tabular*}\caption{Optimal sample fractions, mean values and mean squared errors of the estimator at its optimal levels, for the AR(1) process, with $\eta=0.75,$ $\eta=0.89$ and $\eta=0.96$.}
\end{center}
\end{table}

\subsection{Some overall conclusions}

\begin{itemize}
\item The sample paths of the runs estimator of the upcrossings index have very different patterns for the ARMAX and the AR(1) process. Whilst for the ARMAX process the estimates increase with the value of $k,$ for the AR(1) process the estimates, as function of $k,$ have almost a symmetric distribution, decreasing with smaller values of $k.$
\item  For the ARMAX process, smaller values of $\eta$ tend to be associated with a wider ``bathtub'' pattern of the mean squared error as a function of $k.$
\item  The runs estimator of the upcrossings index has some mean value stability around the target $\eta$ for small values of $k,$ exhibiting, for such values, the mean squared error a ``bathtub'' pattern, although not very wide.
\item For small values of $k,$ or equivalently high levels, we obtain good estimates for $\eta.$ We remark that the choice of the number of top order statistics is a complex problem in extreme value applications.
\item The sensitivity of the mean value to the changes in $k$ seem to clarify the need of studying the bias properties of this estimator. This is clear form Figure 6 where we can see that for the ARMAX process the behaviour of the mean squared error is almost determined by the bias (overlap of both curves), since the variance is always very small. The same behaviour holds for the AR(1) process.
\begin{figure}[!ht]
\begin{center}
\includegraphics[scale=0.3]{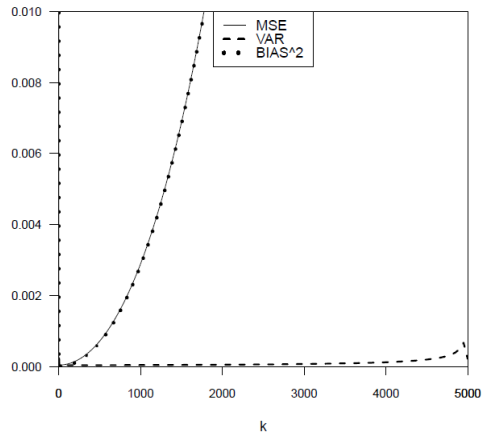}\vspace{-0.5cm}
\caption{Estimated mean squared error (solid line), estimated variance (dashed line) and estimated bias (dotted line), for samples of size $n=5000$ from the ARMAX process in (\ref{armax}).}
\end{center}
\end{figure}

\end{itemize}

\section{Case-studies}\setcounter{equation}{0}

\subsection{Ozone pollution}

\pg We now consider the performance of the above mentioned estimator in the analysis of $n=120$ weekly maxima of hourly averages of ozone concentrations measured in parts per million, in the San Francisco bay area, San José. These data are available in the package Xtremes (Reiss and Thomas \cite{rei}) and have already been studied, for instance, in Gomes {\it{et al.}} \cite{gom} when estimating the extremal index. In Figure 7 we picture the data over the above mentioned period.

\begin{figure}[!ht]
\begin{center}
\includegraphics[scale=0.5]{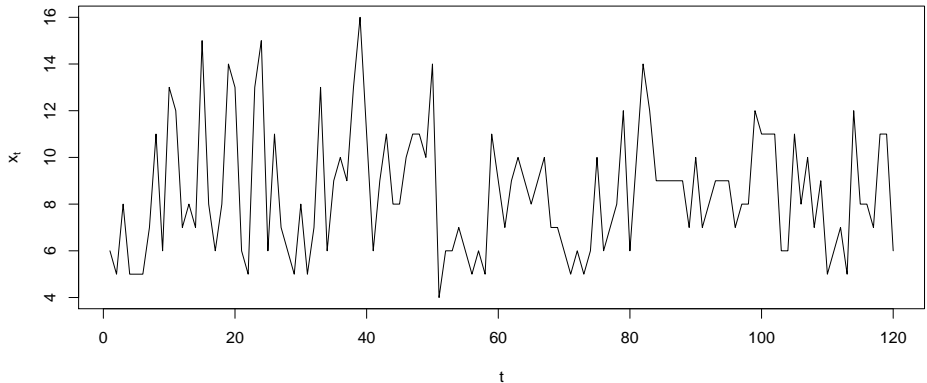}\vspace{-0.5cm}
\caption{Weekly maxima of hourly averages of ozone concentrations measured in parts per million, in the San Francisco bay area, San José}
\end{center}
\end{figure}

Since in practice condition $\widetilde{D}^{(3)}$ is not yet possible to verify we assumed that it holds since as stated Gomes {\it{et al.}} \cite{gom}, most of the parametric models that adequately fit this type of meteorological data satisfy condition $D''$ and therefore also satisfy condition $\widetilde{D}^{(3)}.$

In Figure 8 we picture the sample paths of $\widehat{\eta}_n(k)$ as functions of $k$ in a linear scale (left) and logarithmic (right) scale.

\begin{figure}[!ht]
\begin{center}
\includegraphics[scale=0.25]{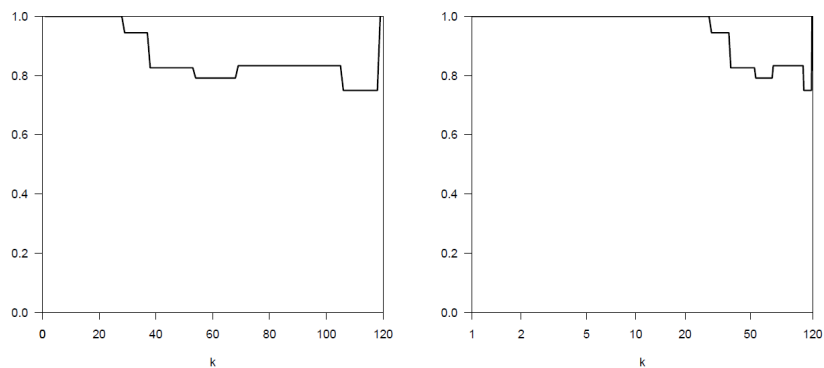}\vspace{-0.5cm}
\caption{Estimates of the upcrossings index as a function of $k$ for the weekly maxima of hourly averages of ozone concentrations measured in parts per million, in the San Francisco bay area, San José, in a linear scale (left) and in a logarithmic scale (right).}
\end{center}
\end{figure}

The stability around one for small values of $k$ agrees with the fact that $\eta=1$ since from \cite{gom} condition $D''$ holds.  Nevertheless, the small number of observations makes it difficult to rely on such a point estimate since the number of upcrossings and consequently the number of non-upcrossings followed by an upcrossing, of a hight level, will always be very small.

\subsection{Financial log-returns}

\pg  Financial time series are very unlikely to satisfy condition $D''$ because of their varying volatility. Therefore, we  shall consider now the time series of log returns of the German stock market index DAX, consisting of the 30 major German companies trading on the Frankfurt Stock Exchange. The daily closing prices of the DAX index for the period from 1991 to 1998 that we work with is available as dataset {\textsf{EuStockMarkets}} in the statistics software R, which we also used for all computations. This data set of 1786 observations has been considered in Klar {\it at al.} \cite{klar} where we can find a full statistical description. In Figure 9 we picture DAX daily closing prices over the mentioned period, $x_t,$ and the log-returns, $100\times (\ln x_t-\ln x_{t-1}),$ the data to be analyzed.

\begin{figure}[!ht]
\begin{center}
\includegraphics[scale=0.25]{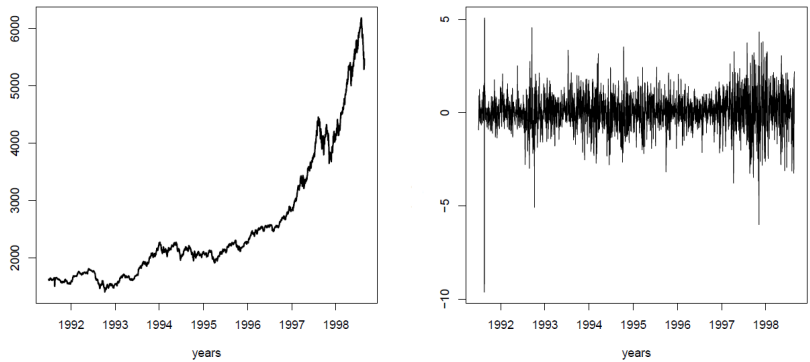}\vspace{-0.5cm}
\caption{DAX daily closing prices (left) and daily log-returns  (right)}
\end{center}
\end{figure}

This time series of log returns is well modeled by a GARCH(1,1), covariance stationary, process (\cite{klar}). Using a goodness-of-fit test  Klar {\it at al.} \cite{klar}, at a 5\% level, do not reject the hypothesis that the innovations of the process follow a  t-distribution with seven degrees of freedom.

\begin{figure}[!ht]
\begin{center}
\includegraphics[scale=0.25]{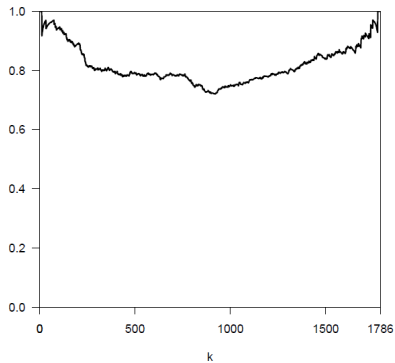}\vspace{-0.5cm}
\caption{Estimates of the upcrossings index as a function of $k$ for the DAX log-returns}
\end{center}
\end{figure}

The graph in Figure 10 allows us to conclude that $\eta$ is not equal to one, which agrees with the fact that these processes are very unlikely to verify condition $D''.$  The sample path exhibits a stability region around a value close to $\eta=0.8.$

Mikosch and St\u{a}ric\u{a} \cite{mik} derive the extremal index for the squared GARCH(1,1) processes, whereas Laurini and Tawn \cite{lau} propose an algorithm for the evaluation of the extremal index of GARCH(1,1) processes with $t-$distributed innovations. It would be interesting, in future work, to obtain similar results for the upcrossings index of a GARCH(1,1) process.

\section{Conclusions}

\pg The upcrossings index, as a measure of the clustering of upcrossings of high levels, is an important parameter when studying extreme events. For sequences satisfying condition $\widetilde{D}^{(3)},$ that locally restricts the dependence of the sequence but still allows clustering of
upcrossings,  we have proposed a simple estimator for this parameter. The study of the properties of the proposed estimator, namely its  consistency and asymptotic normality, has been carried out in Sections 3 and 4. With simulations of well known autoregressive processes that verify condition $\widetilde{D}^{(3)}$ we were able to illustrate the performance of the estimator. Case studies in the fields of environment and
finance were also exploited.

Relation (\ref{rel_eta_teta}) allows one also to estimate $\eta$ through the extremal index $\theta$ modified by consistent estimators of exceedances, $\tau$, and the mean number of upcrossings, $\nu,$ of high levels. Several estimators for the extremal index can be found in the literature (see Ancona-Navarrete and Tawn \cite{anc} for a survey). Other estimators arise from Proposition 1.1 and from the relation with $\mu_1,$ the upcrossings-tail dependence coefficient. In future work we intend to propose new estimators for the upcrossings index and compare the several estimating methods.

%\section*{Agradecimentos}
%We are grateful to the anonymous referee for his detailed comments and suggestions which helped considerably the final form of this paper.

%This research was supported by the research unit ``Centro de Matem\'atica''
%of the University of Beira Interior and the research project PTDC/MAT/108575/2008 through the
%Foundation for Science and Technology (FCT) and co-financed by FEDER/COMPETE.

%% --------------------------------------------------------------
% Refer{\^e}ncias
% Lista das referencias por ordem alfab{\'e}tica de autores-ano (Siga por favor o
% formato apresentado no exemplo):
%% --------------------------------------------------------------

\end{document}